\documentclass[12pt,reqno]{amsart}
\usepackage{verbatim}
\usepackage{amscd}
\usepackage{amsfonts,graphicx,epstopdf,caption}

\begin{document}

\numberwithin{equation}{section}
\renewcommand{\theequation}{\thesection.\arabic{equation}}
\setcounter{secnumdepth}{2}
%

\newcommand{\Bc}{{\mathcal B}}
\newcommand{\Gma}{{\Gamma}}
\newcommand{\Hc}{{\mathcal H}}

\newcommand{\Lc}{{\mathcal L}}

\newcommand{\nab}{\nabla}

\newcommand{\Om}{\Omega}
\newcommand{\Omb}{\overline{\Omega}}

\newcommand{\pal}{\partial}

\newcommand{\Rb}{\overline{\R}}

\newcommand{\R}{\mbox{$\mathbb R$}}

\newcommand{\Sc}{{\mathcal S}}
\newcommand{\Sct}{\tilde{\Sc}}

\newcommand{\sg}{\sigma}
\newcommand{\sgt}{\tilde{\sigma}}

\newcommand{\ubar}{\overline{u}}
\newcommand{\util}{\tilde{u}}

\newcommand{\vap}{\varphi}

\newcommand{\barr}{\begin{eqnarray}}
\newcommand{\bc}{\begin{center}}
\newcommand{\beq}{\begin{equation}}
\newcommand{\bpf}{\begin{proof} \quad}
\newcommand{\btm}{\begin{thm}}

\newcommand{\earr}{\end{eqnarray}}
\newcommand{\ec}{\end{center}}
\newcommand{\eeq}{\end{equation}}
\newcommand{\epf}{\end{proof}}
\newcommand{\etm}{\end{thm}}

\newcommand{\deq}{:= }
\newcommand{\deqs}{\ :=\ }

\newcommand{\eqs}{\ =\ }
\newcommand{\geqs}{\ \geq \ }
\newcommand{\leqs}{\ \leq \ }
\newcommand{\mns}{\, - \, }
\newcommand{\pls}{\, + \, }
\newcommand{\plms}{\ + \ }


\newcommand{\foral}{\qquad \mbox{for all} \quad }
\newcommand{\foreach}{\qquad \mbox{for each} \quad }

\newcommand{\wrt}{ \mbox{with respect to}  }
\newcommand{\xand}{\quad \mbox{and} \quad }
\newcommand{\xfor}{ \quad \mbox{for} \ }
\newcommand{\xiff}{\ \mbox{if and only if} \ }

\newcommand{\xon}{\qquad \mbox{on} \quad}
\newcommand{\xor}{\qquad \mbox{or} \quad}
\newcommand{\xthen}{\quad \mbox{then} \quad}
\newcommand{\xwhen}{\qquad \mbox{when} \quad}
\newcommand{\xwith}{\qquad \mbox{with} \quad}

\newcommand{\xHSEF}{\mbox{harmonic Steklov eigenfunction} }
\newcommand{\xSvl}{\ \mbox{Steklov eigenvalue}}

\newcommand{\ang}[1]{\langle#1\rangle}

\newcommand{\bdy}{\partial \Omega}
\newcommand{\bdRh}{\pal R_h}
\newcommand{\Cinty}{C^{\infty} }

\newcommand{\delj}{\delta_j}
\newcommand{\deljh}{\tilde{\delj}}
\newcommand{\delk}{\delta_k}
\newcommand{\delo}{\delta_1}
\newcommand{\Dbj}{\Dnu b_j}
\newcommand{\Dej}{ \Dnu e_j    }
\newcommand{\Dnu}{D_{\nu}}
\newcommand{\dsg}{\, d \sigma}
\newcommand{\dsgt}{\; d \sgt}
\newcommand{\dxy}{\, dx dy}

\newcommand{\Ebg}{E_b \, g}
 \newcommand{\gamu}{\gamma(u)}
  \newcommand{\gamv}{\gamma(v)}
 \newcommand{\gbar}{\overline{g}}
 
 \newcommand{\gjh}{\hat{g}_j}
 
 \newcommand{\Gamo}{\Gma_1}
 \newcommand{\Gamtw}{\Gma_2}
  \newcommand{\Gamth}{\Gma_3}
  \newcommand{\Gamf}{\Gma_4}

\newcommand{\gradu}{\nabla u}
\newcommand{\gradv}{\nabla v}

\newcommand{\Iby}{\int_{\bdy} \; }
\newcommand{\IMby}{|\bdy|^{-1} \, \Iby \ }
\newcommand{\IOm}{\int_{\Om} \; }

\newcommand{\Lap}{\Delta \, }
\newcommand{\Laps}{\Delta \, s}

\newcommand{\Lapu}{\Delta \, u}
\newcommand{\Lapv}{\Delta \, v}
\newcommand{\Lapw}{\Delta \, w}

\newcommand{\Omo}{\Om_1} 
\newcommand{\OmL}{\Om_L}      

\newcommand{\opal}{\oplus_{\pal}}

\newcommand{\ra}{\rightarrow}

\newcommand{\RN}{{\R}^N}
\newcommand{\Rte}{\R^3}
\newcommand{\Rto}{\R^2}

\newcommand{\ScM}{\Sc_M}

\newcommand{\sj}{s_j}

\newcommand{\sjt}{\tilde{\sj}}
\newcommand{\sk}{s_k}

\newcommand{\sumM}{\sum_{j=1}^M}
\newcommand{\sumo}{\sum_{j=1}^{\infty}}
\newcommand{\sumz}{\sum_{j=0}^{\infty}}

\newcommand{\ujh}{\hat{u}_j}


\newcommand{\COmb}{C(\Omb)}

\newcommand{\Cc}{C_c^1 (\Om)}
\newcommand{\Ccinty}{C_c^{\infty} (\Om)}

\newcommand{\Harm}{\Hc (\Om)}
\newcommand{\Harmh}{{\Hc}^{1/2}(\Om)}
\newcommand{\Honeharm}{{\Hc}^{s+1/2}(\Om)}

\newcommand{\Hone}{H^1(\Om)}
\newcommand{\Wop}{W^{1,p}(\Om)}

\newcommand{\Hzone}{H_{0}^1(\Om)}

\newcommand{\Lone}{L^1(\Om)}
\newcommand{\Lp}{L^p (\Om)}
\newcommand{\Lt}{L^2 (\Om)}

\newcommand{\Cbdy}{C(\bdy)}
\newcommand{\CRh}{C(\pal R_h)}

\newcommand{\Hhby}{H^{1/2}(\bdy)}
\newcommand{\Hmsby}{H^{-s}(\bdy)}
\newcommand{\Hsby}{H^s(\bdy)}

\newcommand{\Lpb}{L^p (\bdy, \dsg)}

\newcommand{\Ltby}{L^2 (\bdy, \dsg)}
\newcommand{\Lttby}{L^2 (\bdy, \dsgt)}

%
\newtheorem{thm}{Theorem}[section]
\newtheorem{cor}[thm]{Corollary}
\newtheorem{cond}{Condition}
\newtheorem{lem}[thm]{Lemma}
\newtheorem{prop}[thm]{Proposition}

\title[ Steklov Approximations of Harmonic Functions]
{Steklov Approximations of Harmonic  Boundary Value Problems on Planar Regions.}
\author[Auchmuty and Cho]{Giles Auchmuty and Manki Cho}

\address{Department of Mathematics,  University of Houston, Houston, Tx 77204-3008}
\email{auchmuty@uh.edu}
\address{School of Mathematical Sciences, Rochester Institute of Technology, Rochester, New York, 14623}
\email{mxcsma1@rit.edu}

\thanks{The first author gratefully acknowledges research support by NSF award DMS 11008754. 
 }
\date{\today; \quad SAHFin.tex}

\begin{abstract} 
Error estimates for approximations of harmonic functions on planar regions by subspaces spanned by 
the first $\xHSEF$s  are found.  
They are based on the explicit representation of harmonic functions in terms of these $\xHSEF$s.
When the region is a rectangle of aspect ratio h, some computational results regarding these approximations
for problems with known explicit solutions are described.
\end{abstract}
\maketitle


\section{Introduction}\label{s1}

This paper will describe results about the approximation of solutions of various boundary value problems for 
Laplace's equations on a bounded planar regions  $\Om$ in the plane  using certain harmonic Steklov eigenfunctions.
That is our interest is in approximating harmonic functions $u: \Omb \to \R$ that satisfy either Dirichlet, Robin or 
Neumann boundary conditions 
\beq \label{BCs}
u \eqs g \xor \Dnu u \pls b \, u \eqs g \xon \bdy \eeq 
Here $\nu$ is the outward unit normal and $b \geq 0$ is a constant. 

It has been shown that there are orthogonal  bases of the class of all finite energy harmonic  functions $\Harm$ on $\Om$
consisting of harmonic Steklov eigenfunctions - as summarized below in section \ref{s3}.
A natural question is how good are  approximations using specific subclasses of these eigenfunctions; particularly those
corresponding to the lowest Steklov eigenvalues?
Here some general results about these approximations are described in sections \ref{s4} and \ref{s6} and some
computational results for particular problems with exact solutions are described in sections \ref{s5} and \ref{s7}.
The computational results are based on the fact that explicit formulae are known for the Steklov eigenvalues and eigenfunctions
on rectangles of arbitrary aspect ratio h.

Existence-uniqueness theorems for these problems may be found in most texts that treat elliptic boundary value problems.
Under differing assumptions on $g, \Om$ and $\bdy$ the solutions $u$ are $\Cinty$ on $\Om$ and lie in various different
Banach or Hilbert spaces of functions on $\Om$.
For an excellent review of classical results about  these problems see chapter 2 of \cite{DL} by Benilan. 
In particular a function $u \in \Lone$ is said to be an {\it ultraweak} solution of Laplace's equation provided it obeys
\beq \label{harm}
\IOm \ u \ \Lap \vap \ \dxy \eqs 0 \foral \vap \in \Ccinty. \eeq
Such an ultraweak solution  is a {\it classical solution} of  Laplace's equation provided it  is equivalent to a continuous 
function on $\Omb$. 
Such solutions need not be weak solutions in the Sobolev space $\Hone$ - even when $\Om$ is a disk in the plane.

In the  sections \ref{s3},  \ref{s4} and \ref{s6},  general results about Steklov approximations of harmonic functions 
 will be described.
These are based on the observation that there is an algorithm for constructing a basis of classes of harmonic 
functions on regions consisting of  harmonic Steklov eigenfunctions. 
See Auchmuty \cite{Au1}-\cite{Au3} and the boundary traces of  these eigenfunctions are bases of associated Hilbert
 spaces of functions on $\bdy$.
In particular various error estimates for Steklov approximations are obtained. 

When the domain is a planar disk, the Steklov eigenfunctions are the usual harmonic functions $r^m \cos{m \theta}, r^m
\sin{m \theta}$ of Fourier analysis and the question of the approximation of harmonic functions on the unit  disc by
harmonic polynomials has a huge literature.
The text of Axler, Bourdon and Ramey \cite{ABR} is a recent introduction to the theory. 

The use of Steklov bases for spaces of harmonic functions permits the generalization of some the results of 
classical harmonic function theory to quite general bounded (and even exterior) regions in $\RN$, see \cite{Au3}.
Here some results about the approximation of harmonic functions by finite sums of Steklov eigenfunctions will be
described with computational results for the case of a rectangle.

When the region $\Om$ is a rectangle, the Steklov eigenfunctions are known explicitly see Auchmuty and Cho \cite{AC}
or Girouard and Polterovich \cite{GP} where a completeness proof may be found.
The paper \cite{AC} described the generalization of the mean value theorem to rectangles and to cases where Robin data 
on $\bdy$ is known.
Here  in sections \ref{s5} and \ref{s7} computational results for  Steklov approximations of certain harmonic functions
regarded as solutions of Laplace's equations with  various boundary value conditions are described.

For general regions, the Steklov eigenvalues and eigenfunctions are not (yet) known explicitly.
However a number of authors have studied the numerical determination of these eigenfunctions including
Cheng, Lin and Zhang \cite{CLZ}, and Kloucek, Sorensen and  Wightman \cite{KSW}. 
The software FreeFem++ \cite{He}  has subroutines for the computation of Steklov eigenfunctions and eigenvalues that was
used for confirmation of some of the  analytical results described here.

After introducing our assumptions and notation in section 2, the basic properties of harmonic Steklov eigenproblems are
described in section 3.
In particular the formula for the ($L^2-$)harmonic extension operator - or solution operator for the Dirichlet harmonic  problem
is described. 
In section 4 some error estimates for Steklov approximations are given.
Then some computational results about such  problems on a rectangle are described in section 5.
 Dirichlet problems are considered  in sections 4  and 5 and then  results for Robin and Neumann problems are described
  in sections 6 and 7.
  
  Our general conclusion is that the use of Steklov approximations of harmonic functions should be seriously considered 
  by researchers needing good approximations using a small dimensional subspace of harmonic functions.
  Once the Steklov eigenfunctions ane eigenvalues of a region are known, simple formulae provide excellent approximations of the solutions in a region.
   Properties of the approximations  close to, or on, the boundary  needs significant further study.
   It should be noted that this analysis  extends to the solution much more general self-adjoint second order elliptic equations
   of the form $\Lc u = 0$  using similar  general constructions.

\vspace{1em}
%

\section{Assumptions and Notation.}\label{s2}

This paper treats various Laplacian boundary value problems on regions  $\Om$ in the plane $\Rto$.
A region  is a non-empty, connected,  open subset of $\Rto$. 
Its closure  is denoted  $\Omb$ and its boundary is  $\bdy \deq \Omb\setminus \Om$. 
Some regularity of the boundary $\bdy$ is required.
Each component (= maximal connected closed subset) of the boundary is assumed to be   a Lipschitz continuous closed curve.
Let $\sigma$ denote arc-length along a curve so the  unit outward normal $\nu(z)$  is defined $\sigma \, a.e.$

 $\Lp$ and $\Lpb, \ 1 \le p \le \infty$ are the usual spaces with  p-norm denoted by
${\|u\|_p}$ or $ \|u \|_{p,\bdy}$ respectively.
 When $p=2$ these are real Hilbert spaces with inner products defined by 
\[ \ang{u , v}  \deq \int_{\Om} \  u(x) \; v(x) \ dx \qquad \mbox{and} \quad
 \ang{u , v}_{\bdy}  \deqs |\bdy|^{-1} \ \int_{\bdy} \  u \; v \ \dsg.  \]
 $\COmb$ is the space of continuous functions on the closure $\Omb$ of $\Om$ with the sup norm
 $\| u \|_b := sup_{\Omb} \ |u(x,y)| $. 

The weak j-th  derivative of u is  $D_j u$ - and all derivatives will be taken in a weak sense.
 Then $\gradu(x) := \, (D_1u(x), \ldots, D_Nu(x))$ is the gradient of u  and $\Hone$ is the usual real Sobolev space 
 of  functions  on $\Om$. 
 It is a real Hilbert space under  the standard $H^1-$ inner product
\begin{equation}\label{ip1}
{[u , v]}_1  \deq \IOm \ [ u(x) . v(x) \  + \  \gradu(x)  \cdot  \gradv(x)] \  dx.
\end{equation}
The corresponding norm   is denoted   $ \|u \|_{1,2}$.

The region $\Om$ is said to satisfy {\it Rellich's theorem} provided  the imbedding of   $\Hone$ into $\Lp$
 is compact for $1 \leq p < \infty$.

The boundary trace operator  $\gamma:  \Hone \to \Ltby$ is the linear extension of the map restricting  Lipschitz 
continuous functions on $\Omb$ to $\bdy$. 
The region $\Om$ is said to satisfy a {\it compact trace theorem} provided  the  boundary trace mapping 
$\gamma : \Hone \ra \Ltby$ is compact. 
Usually $\gamma$ is omitted so $u$ is used in place of $\gamu$ for the  trace of a function on $\bdy$.

The {\it Gauss-Green} theorem   holds on $\Om$ provided 
\beq \label{GG} 
\IOm \  u(x) D_j v(x) \  dx \eqs  \int_{\bdy} \ \gamu \; \gamv \, \nu_j \,  \dsg \ - \ \IOm \  v(x) \, D_j  u(x) \,  dx \xfor \ 1 \le j \le N.
\eeq
for all $u,\ v$ in $\Hone$. 
The  requirements on the region will be 

\noindent{\bf Condition B1:}  \quad {\it $\Om$ is a bounded region in  $\Rto$  whose boundary $\bdy$ is  a finite number of
disjoint closed Lipschitz curves, each of finite length  and such that   the Gauss-Green,  Rellich and compact trace theorems hold.}

We will use the   equivalent inner products on $\Hone$ defined by
\begin{equation}\label{ippal}
{[u , v]}_{\pal}  \deq \IOm \  \gradu  \cdot \gradv \  dx \pls  \int_{\bdy} \  u \; v \ \dsg.
\end{equation}
The corresponding norm will be denoted by  $ \| u \|_{\pal}$.
The proof that this  norm is equivalent to the usual $(1,2)-$norm on $\Hone$ when  (B1) holds is Corollary 6.2 of \cite{Au1} and also 
is part of theorem 21A of \cite{Ze}. 

A function $u \in \COmb$ or $\Hone$ is said to be harmonic provided it satisfies \eqref{harm}.
Define $\Harm$ to be the space of all  harmonic functions in $\Hone$. 
When (B1)  holds, the closure of $\Cc$ in the $H^1-$norm is the usual Sobolev space $\Hzone$.
Then (\ref{harm}) is equivalent to saying  that $\Harm$ is $\pal-$orthogonal to $\Hzone$. 
This may be expressed as
 \beq \label{e2.5}
\Hone \eqs \Hzone \oplus_{\pal}  \Harm,
\eeq
where $\oplus_{\pal}$ indicates that this is a $\pal-$orthogonal decomposition.
 
The analysis to be described here is based on the construction of a $\pal-$orthogonal basis of the Hilbert space
$\Harm$ consisting of harmonic Steklov eigenfunctions. 
In particular we shall prove results about the approximation of solutions of harmonic boundary value problems by such
eigenfunctions.

\vspace{1em}
\section{Steklov Representations of Solutions of Harmonic Boundary Value Problems. } \label{s3}

Let $\Om$ be a bounded region in $\Rto$ that satisfies (B1).
 A non-zero function $s \in \Hone$  is said to be a {\it harmonic Steklov eigenfunction}  on $\Om$ 
 corresponding to the  Steklov eigenvalue $\delta$ provided $s$ satisfies
\beq \label{HSeqn}
\IOm \;  \nabla s  \cdot  \nabla v  \; dx \eqs \delta  \   \ang{s,v}_{\bdy} \eqs \delta \ \IMby s \, v \dsg.  \qquad 
\mbox{for all} \quad v \in \Hone.
\eeq

  This is the weak form of the boundary value problem
\beq \label{Weeqn}
\Delta \, s \eqs 0 \qquad \mbox{on $\Om$ with} \quad  \Dnu \, s \eqs \delta \ |\bdy|^{-1} \  s 
\quad \mbox{on} \ \bdy.
\eeq
Here $\Delta$ is the Laplacian and $\Dnu \, s(x) := \nabla s(x) \cdot \nu(x)$ is the unit 
outward normal derivative of $s$    at a point on the boundary. 

Descriptions of the analysis of these eigenproblems may be found in Auchmuty \cite{Au1} - \cite{Au4}.
These eigenvalues and a corresponding family of $\pal-$orthonormal eigenfunctions 
may be found using variational principles as described in sections 6 and 7 of Auchmuty \cite{Au1}. 
$\delta_0 = 0$ is the least eigenvalue of this problem corresponding to the eigenfunction
$s_0 (x) \equiv 1$ on $\Om$. This eigenvalue is simple as $\Om$ is connected. 
Let the first k Steklov eigenvalues be $0 = \delta_0 < \delta_1 \le \delta_2  \leq \ldots \leq \delta_{k-1}$ 
and $s_0, s_1, \ldots, s_{k-1}$ be a corresponding  set of  $\pal-$orthonormal eigenfunctions. 
The k-th eigenfunction $s_k$ will be a maximizer of the functional
\beq \label{e3.3}
\Bc(u) \deqs  \Iby \| \gamma( u) |^2 \, \dsg
\eeq
over the subset $B_k$ of functions in $\Hone$ which satisfy 
\beq \label{e3.5}
\|u\|_{\pal} \leqs 1 \qquad \mbox{and} \qquad \  \ang{\gamma( u), \gamma (s_l) \, }_{\bdy} \eqs 0 \quad 
\mbox{for} \quad 0  \leq l \leq k-1.
\eeq

The existence and some properties of such eigenfunctions are described in sections  6 and 7 of \cite{Au1} 
for a more general system. 
In particular, that analysis shows that  each $\delta_j$ is of finite multiplicity and $\delta_j \ra \infty$ as $j \ra \infty$; 
see Theorem 7.2 of   \cite{Au1}.   
The maximizers not only are $\pal -$orthonormal but they also satisfy
\barr 
\IOm   \nabla s_k \cdot  \nabla s_l  \; dx \eqs  \IMby \,  s_k \, s_l  \, \dsg \eqs 0 \qquad 
\mbox{for } \quad k \ne l. \label{e3.7} \\
\IOm  | \nabla s_k |^2  \; dx \eqs \frac{\delk}{1+\delk} \qquad \text{and} \quad 
 \IMby \,  |\gamma (s_k) |^2   \, \dsg \eqs \frac{1}{1+\delk} \quad \mbox{for} \quad k \ge 0.  \label{e3.8}
\earr
Recently  Daners \cite{Da}  corollary 4.3 has  shown that, when $\Om$ is a Lipschitz domain, then the Steklov eigenfunctions are 
continuous on $\Omb$. 

The analysis in this paper is based on the fact that  harmonic Steklov eigenfunctions on $\Om$ can be chosen 
to be  orthogonal  bases of both  $\Harm$ and of $\Ltby$.  
It should be noted that, for regions other than discs (or balls in higher dimenions),  these Steklov eigenfunctions are generally 
 not $L^2-$orthogonal on $\Om$. 

Let $\Sc := \{s_j : j \geq 0 \}$ be the maximal family of  $\pal -$orthonormal  eigenfunctions constructed inductively as above. 
For this paper,  it is more convenient to use the Steklov eigenfunctions normalized by their boundary norms. \\
Define the functions $\sjt(x) \deqs \sqrt{1+\delj} \quad  \sj(x)$ for $j \geq 0$. 
From \eqref{e3.8}, these satisfy
\beq \label{e3.10}
 \Iby \tilde{s_j} \, \tilde{s}_k \dsg \eqs 0 \xwhen j \neq k \xand  \Iby {\tilde{s_j}}^2  \dsg \eqs |\bdy|. 
   \eeq
 These Steklov eigenfunctions are said to be {\it boundary normalized} and the associated set 
 $\Sct \deqs \{\sjt : j \geq 0\} $ is an orthonormal basis of $\Ltby$. See theorem 4.1 of \cite{Au2}.

For given $g \in \Ltby$, let 
\beq \label{e3.11}
g_M (x,y) \deqs \gbar \pls \sum_{j=1}^M \ \gjh \ \sjt(x,y) \xwith \gjh = \ang{g, \sjt}_{\bdy} \eeq
 be the M-th Steklov approximation of $g$ on $\bdy$. 
Here $\gbar := g_0$ is the mean value of $g$ on $\bdy$ and this is the standard projection  of elements in a Hilbert space 
onto subsets of an  orthonormal basis
Note that $g_M$ is continuous and bounded  on $\bdy$ as each $\sjt$ is and $g_M$ converges strongly to $g$ in $\Ltby$
from the Riesz-Fischer theorem and 
\beq \label{e3.12}
\| \, g - g_M \, \|_{2,\bdy}^2 \eqs \| \, g \, \|_{2,\bdy}^2 \mns  \| \,  g_M \, \|_{2,\bdy}^2. \eeq

The unique solution of Laplace's equation on $\Om$ subject to the Dirichlet boundary condition $\gamu = g$ on $\bdy$ is given by
\beq \label{e3.15}
u(x,y) \eqs E_H g(x,y) \eqs \gbar \pls \lim_{M \to \infty} \  \sum_{j=1}^M \ \gjh \sjt(x,y) \xfor (x,y) \in \Om. \eeq
See section 6 of \cite{Au3} for a proof. $E_H$ will be called the {\it harmonic extension} operator and is
a compact linear map from $\Ltby$ to $\Lt$.
Classically this map has been represented as an integral operator with  the {\it Poisson kernel}.
Theorem 6.3 of \cite{Au3} says that $E_H$ is an isometric isomorphism of $\Ltby$ with a space denoted $\Harmh$ that is a proper subspace of $\Lt$.

\vspace{2em}
\section{Error Estimates for Steklov Approximations } \label{s4}

First some error estimates for the Steklov approximations of solutions of the Dirichlet problem for Laplace's equation may be
described. 
Essentially bounds on the error from these formulae for Steklov approximations are found that depend only on the errors
in the Steklov approximations of $g$ on the boundary.
Let
\beq \label{e4.1}
u_M (x,y) \deqs \ubar \pls   \sum_{j=1}^M \ c_j \ \sj(x,y) \xfor (x,y) \in \Om. \eeq
be a finite sum of the first M+1 $\xHSEF$s on $\Om$. 
Then for each integer $M, \ u_M$ is $\Cinty$ on $\Om$, continuous on $\Omb$ and in $\Harm$.
The subspace spanned by $\ScM := \{\sj : 0 \leq j \leq M\} $ will be denoted $V_M$. 

\begin{thm} \label{T4.1}
Assume  $\Om, \bdy$ satisfy (B1) and $\Sc, \Sct$ are the orthonormal bases of $\Harm$, $\Ltby$ described above.
If $g \in \Cbdy, u = E_H g$ and $ u_M$ is defined by  \eqref{e4.1}, then
\beq \label{e4.3} 
\| u \mns u_M \, \|_{2,\Om} \leqs C_2 \ \| g \mns g_M \, \|_{2,\bdy} \xand \| u \mns u_M \, \|_{\infty,\Om} \leqs \| g \mns g_M \, \|_{\infty,\bdy} 
\eeq
where $C_2$ is the Fichera constant for $\Om$ and $g_M = \gamma(u_M)$.  \etm
\bpf
Here $g_M$ is the boundary trace of $u_M$, so $u_M$ is the harmonic extension of $g_M$.
The 2-norm inequality is Fichera's inequality with $C_2$ being the first eigenvalue of the Dirichlet, biharmonic Steklov eigenproblem.
See \cite{Fic} for the original version and \cite{Au5} for a recent description and proof under weak boundary regularity
conditions. 
The second inequality is the maximum principle for classical solutions of Laplace's equation.  
\epf

It should be noted here that the inequalities in \eqref{e4.3} do not require that $g_M$ be the Steklov approximation of g on $\bdy$. 
They hold for any function that is a linear combination of the first {M+1}  Steklov eigenfunctions. 
Note that $L^p-$bounds follow for $2 < p < \infty$ by interpolation.
When $g_M$ is the M-th Steklov approximation of $g$ on $\bdy$ as in \eqref{e3.11}, then one also has
 
\btm \label{T4.2}
 Assume (B1) and $g \in \Hhby, \ g_M$  is defined by \eqref{e3.11}, $u = E_H g$ and $u_M = E_H g_M$. 
 Then $g_M$ converges strongly to $g$ in $\Hhby$ and $u_M$ converges uniformly to $u$ on compact subsets of $\Om$.
 Moreover 
 \beq \label{e4.5}
 \| \nabla (u - u_M) \, \|_{2,\Om}^2 \eqs \sum_{j = M+1}^{\infty} \ \delj \gjh^2 \eqs \|g \|_{1/2, \bdy}^2  \mns \|g_M \|_{1/2, \bdy}^2
 \eeq   \etm
\bpf
The fact that $g_M$ converges strongly to $g$ in $\Hhby$ and $\Hone$ follows from the fact that $\Sc$ is an
orthonormal basis of $\Harm$. 
The proof of uniform convergence is standard, while \eqref{e4.5} follows from the orthogonality properties of 
Steklov eigenfunctions. 
\epf

Also note that the Steklov eigenfunction have scaling properties. 
 Given  $\Omo \subset \Rto$, let $\OmL := \{ Lx : x \in \Omo \}$ with $L > 0$. 
When  $h$ is a harmonic function on $\Omo$, then the function $h_L(y) := h(y/L) $ will be a harmonic function on $\OmL$.
If $h$ is a harmonic Steklov eigenfunction on $\Omo$  with  $\xSvl \ \delta$, then $h_L$ will be a $\xHSEF$  on 
$\OmL$ with the $\xSvl \ \delta / L $.
Thus it  suffices to study problems with a normalized bounded region $\Omo$; the eigenvalues and eigenfunctions for scalings of a
region then follow from these formulae.

The following sections will look at some aspects of the approximation of solutions of  Laplace's equation on
rectangles by finite sums of the form \eqref{e3.15}.
Rectangles are chosen since  we have explicit expressions for the Steklov eigenfunctions and eigenvalues on rectangles.

\vspace{2em}
\section{Steklov Approximations of Harmonic Functions on a Rectangle } \label{s5}

When $\Om \eqs R_h \deqs (-1,1) \times (-h,h)$ is a rectangle with aspect ratio $h$, the Steklov eigenfunctions and eigenvalues are known explicitly. See Auchmuty and Cho \cite{AC} section 4 where eight families of eigenfunctions are described  and characterized by their symmetry properties with respect to the center. Class I eigenfunctions are even in x and y, class II are odd in x and y, class III are even in x and odd in y, class IV are odd in x and even in y.

\begin{figure}[here!]
\begin{center}
 \includegraphics[height=6.5cm, width=12cm]{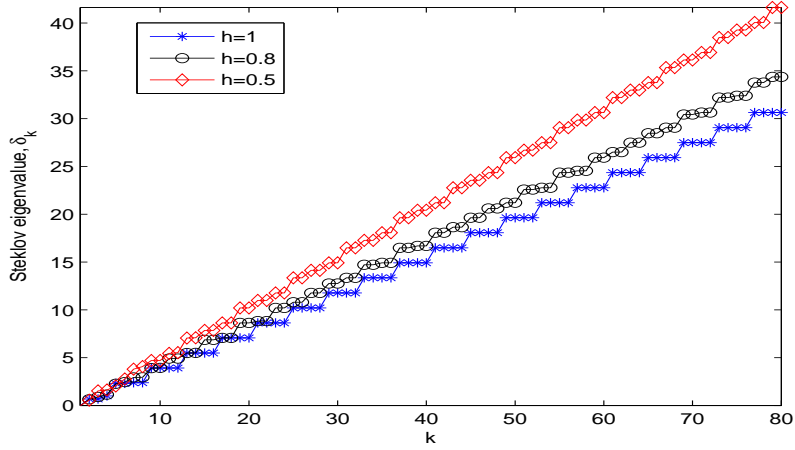}
 \label{eig}
  \caption{First 80 Steklov eigenvalues on $R_h$ corresponding to $h=1, 0.8,$ and $0.5$}
  \end{center}
\end{figure}

By separation of variables the explicit formulae for the Steklov eigenfunctions may be found.
The first  eigenfunction $s_0(x,y) \equiv 1$ is in class I and the other (unnormalized) Steklov eigenfunctions have the forms
\begin{eqnarray} 
s(x,y) \deqs \cosh{\nu x} \, \cos{\nu y} &\xwhen& \quad \tan{\nu h} \pls \tanh{\nu} \eqs 0, \label{s4.1} \\
s(x,y) \deqs \cos{\nu x} \, \cosh{\nu y} &\xwhen& \quad \tan{\nu} \pls \tanh{\nu h} \eqs 0. \label{s4.2}
\end{eqnarray}

When $h =1$, the first eigenfunction in class II is $s_3(x,y) = xy$. 
Otherwise the (unnormalized)  eigenfunctions and eigenvalues in this class have the forms 
\begin{eqnarray} 
s(x,y) \deqs \sinh{\nu x} \, \sin{\nu y} &\xwhen& \quad \cot{\nu h} \mns  \coth{\nu} \eqs 0,  \label{s4.3}\\
s(x,y) \deqs \sin{\nu x} \, \sinh{\nu y} &\xwhen& \quad \cot{\nu} \mns  \coth{ \nu h}\eqs 0 . \label{s4.4}
\end{eqnarray}

Similarly eigenfunctions in class III have the forms
\begin{eqnarray} 
s(x,y) \deqs \cosh{\nu x} \, \sin{\nu y} &\xwhen& \quad \cot{\nu h} \mns \tanh{\nu} \eqs  0, \label{s4.5} \\
s(x,y) \deqs \cos{\nu x} \, \sinh{\nu y} &\xwhen& \quad \tan{\nu} \pls \coth{\nu h} \eqs 0 \label{s4.6}
\end{eqnarray}

Finally the eigenfunctions in class IV have the forms
\begin{eqnarray} 
s(x,y) \deqs \sinh{\nu x} \, \cos{\nu y} &\xwhen& \quad \tan{\nu h} \pls \coth{\nu} \eqs 0 \label{s4.7} \\
s(x,y) \deqs \sin{\nu x} \, \cosh{\nu y} &\xwhen& \quad \cot{\nu} \mns \tanh{\nu h} \eqs 0. \label{s4.8}
\end{eqnarray}

The associated Steklov eigenvalues, $\delta$ are 
\begin{itemize}
\item[(i)] $\delta =  \nu \tanh{\nu} \quad $ when $\nu$ is a solution of the equation in  \eqref{s4.1} or \eqref{s4.5}.
\item[(ii)] $\delta =  \nu \tanh{\nu h} \quad$ when $\nu$ is a solution of the equation in  \eqref{s4.2} or \eqref{s4.8}.
\item[(iii)] $\delta =  \nu \coth{\nu} \quad $ when $\nu$ is a solution of the equation in   \eqref{s4.3} or \eqref{s4.7}.
\item[(iv)] $\delta =  \nu \coth{\nu h} \quad $ when $\nu$ is a solution of the equation in  \eqref{s4.4} or \eqref{s4.6}. 
\end{itemize}

Knowing these explicit formulae for the eigenvalues and eigenfunctions  the approximations  of some given harmonic 
functions using relatively few harmonic Steklov eigenfunctions will be computed.
Since there are eight families of harmonic Steklov eigenfunctions associated with different even/odd symmetries
about the center we have concentrated on approximations involving the first 8M  eigenfunctions with M = 2, 3 and 5.

Note that the convergence results for the Steklov series expansions hold only when the coefficients are precisely the 
Steklov coefficients $\gjh$ defined by \eqref{e3.11}.
The value of $\ubar$ is the mean value of the integral of $g$ around $\bdy$.
However the approximation results of section \ref{s4} hold quite generally for any choice of coefficients.

For the following calculations the coefficients were obtained by evaluating the boundary integrals $\gjh$ of \eqref{e3.11} using  
the global adaptive quadrature(MATLAB's integral). The absolute and relative error tolerance are $10^{-10}$ and $10^{-6}$, respectively.
Then the M-th Steklov approximation  $u_M$ is the function defined by \eqref{e4.1} with $\ujh = \gjh$.

The following tables illustrate the pointwise approximations obtained for these sums at the points 
$P_1 = (0.9, 0.9), \ P_2 =(0.9,0.1) , \ P_3=(0.8,0.6),\ P_4=(0.3,0.9), \ P_5=(0.5, 0.5)  $ and for M = 2,3, 5 and the exact results to 6 decimal places.
Let $D_M(x,y):=|g(x,y)-g_M(x,y)|$ be the absolute error at $(x,y)$. Also let $f_1(x,y):=x^4-6x^2y^2+y^4$, $f_2(x,y):=\frac{2-x}{(2-x)^2+y^2}$, and $f_3(x,y):=\ln(\sqrt{(x-3)^2+(y-3)^2})$.

\begin{table}[h!]
\begin{center}
\begin{tabular}{|c||c|c|c|c|c| }
\hline
&$P_1$&$P_2$&$P_3$&$P_4$&$P_5$\\ \hline
M=2&$-2.626748$ &$0.694643$& $-0.844238$&$0.230283$& $-0.249859 $\\ \hline
M=3& $-2.625942$& $0.607979$& $-0.842944$& $0.225907$&$ -0.249983$\\ \hline
M=5& $-2.624712$& $0.607588$& $-0.843208$& $0.226837$& $-0.250000$\\ \hline 
$g(x,y)$& $-2.624400$& $0.607600$& $-0.843200$& $0.226800$& $-0.250000$\\ \hline 
$D_2(x,y)$& $ 0.002348$ &$0.002957$&$0.001038$&$ 0.003483$&$ 0.000141$\\ \hline 
$D_3(x,y)$& $ 0.001542$&$ 0.000379$&$ 0.000256$&$ 0.000893 $&$ 0.000017$\\ \hline 
$D_5(x,y)$& $0.000312$&$ 0.000012$&$ 0.000008$&$ 0.000037$&$ 0$\\ \hline 
\end{tabular}
\end{center}
\caption{$g(x,y)=f_1(x,y)$ and $h=1$}
\end{table}

\begin{table}[here!]
\begin{center}
\begin{tabular}{|c||c|c|c|c|c| }
\hline
&$P_1$&$P_2$&$P_3$&$P_4$&$P_5$\\ \hline
M=2&0.544285&0.899505& 0.666815& 0.455438& 0.600096\\ \hline
M=3&0.544745& 0.902138& 0.667202& 0.460368& 0.599985\\ \hline
M=5& 0.544675& 0.901609& 0.666636& 0.459219& 0.600000\\ \hline 
$g(x,y)$&0.544554& 0.901639& 0.666667& 0.459459& 0.600000 \\ \hline 
$D_2(x,y)$& $ 0.000269$ &$0.002135$&$0.000148$&$ 0.004021$&$ 0.000096$\\ \hline 
$D_3(x,y)$& $ 0.000191$&$ 0.000498$&$ 0.000535$&$ 0.000909 $&$ 0.000015$\\ \hline 
$D_5(x,y)$& $0.000121$&$ 0.000030$&$ 0.000031$&$ 0.000240$&$ 0$\\ \hline
\end{tabular}
\end{center}
\caption{$g(x,y)=f_2(x,y)$ and $h=1$}
\end{table}

\begin{table}[h!]
\begin{center}
\begin{tabular}{|c||c|c|c|c|c| }
\hline
&$P_1$&$P_2$&$P_3$&$P_4$&$P_5$\\ \hline
M=2& 1.088867& 1.277069& 1.179619& 1.230746& 1.262756\\ \hline
M=3& 1.088349& 1.274927& 1.180394& 1.229961& 1.262881\\ \hline
M=5& 1.088384& 1.275412& 1.180439& 1.229874& 1.262864\\ \hline 
$g(x,y)$& 1.088511& 1.275503& 1.180427& 1.229794&1.262864 \\ \hline 
$D_2(x,y)$& $  0.000356$& $ 0.001566$& $ 0.000808$& $ 0.000952$& $ 0.000108$\\ \hline 
$D_3(x,y)$& $ 0.000162$& $ 0.000576$& $ 0.000033$& $ 0.000167$& $ 0.000017$\\ \hline 
$D_5(x,y)$& $0.000127$& $ 0.000091$& $ 0.000012$& $ 0.000080$& $ 0$\\ \hline
\end{tabular}
\end{center}
\caption{$g(x,y)=f_3(x,y)$ and $h=1$}
\end{table}

Let rerr$_{\infty}(g):=\frac{||g-g_M||_{\infty,\partial\Omega}}{||g||_{\infty,\partial\Omega}}$ and rerr$_{2}(g):=\frac{||g-g_M||_{2,\partial\Omega}}{||g||_{2,\partial\Omega}}$ be the relative error of $M-$th Steklov approximation of $g$ in $L^{\infty}(\Omega)$ norm and $\Ltby$, respectively.

\begin{table}[here!]
\begin{center}
\begin{tabular}{|c||c|c|c| }
\hline
&rerr$_{\infty}(f_1)$&rerr$_{\infty}(f_2)$&rerr$_{\infty}(f_3)$\\ \hline
M=2&$6.59553\times 10^{-3}$&$1.82382\times 10^{-2}$&$6.48245\times 10^{-3}$\\ \hline
M=3&$2.28748\times 10^{-3}$&$1.21554\times 10^{-2}$&$4.3219\times 10^{-3}$\\ \hline
M=5&$5.55757\times 10^{-4}$&$7.35222\times 10^{-3}$&$2.59338\times 10^{-3}$\\ \hline 
\end{tabular}
\end{center}
\caption{Relative errors of the Steklov approximations of $f_1, f_2,$ and $f_3$, respectively where $h=1$}
\end{table}

\begin{table}[here!]
\begin{center}
\begin{tabular}{|c||c|c|c| }
\hline
&rerr$_{\infty}(f_1)$&rerr$_{\infty}(f_2)$&rerr$_{\infty}(f_3)$\\ \hline
M=2&$4.82556\times 10^{-2}$&$2.46749\times 10^{-2}$&$6.38229\times 10^{-3}$\\ \hline
M=3&$4.20662\times 10^{-2}$&$1.78505\times 10^{-2}$&$4.18945\times 10^{-3}$\\ \hline
M=5&$2.28023\times 10^{-2}$&$1.0105\times 10^{-2}$&$2.47618\times 10^{-3}$\\ \hline 
\end{tabular}
\end{center}
\caption{Relative errors of the Steklov approximations of $f_1, f_2,$ and $f_3$, respectively where $h=0.8$}
\end{table}

\begin{table}[here!]
\begin{center}
\begin{tabular}{|c||c|c|c| }
\hline
&rerr$_{\infty}(f_1)$&rerr$_{\infty}(f_2)$&rerr$_{\infty}(f_3)$\\ \hline
M=2&$2.09505\times 10^{-1}$&$3.40908\times 10^{-2}$&$5.58445\times 10^{-3}$\\ \hline
M=3&$1.12233\times 10^{-1}$&$2.00031\times 10^{-2}$&$3.84456\times 10^{-3}$\\ \hline
M=5&$7.66842\times 10^{-2}$&$1.29479\times 10^{-2}$&$2.24773\times 10^{-3}$\\ \hline 
\end{tabular}
\end{center}
\caption{Relative errors of the Steklov approximations of $f_1, f_2,$ and $f_3$, respectively where $h=0.5$}
\end{table}

\begin{table}[here!]
\begin{center}
\begin{tabular}{|c||c|c|c| }
\hline
&rerr$_{2}(f_1)$&rerr$_{2}(f_2)$&rerr$_{2}(f_3)$\\ \hline
M=2&$5.22051\times 10^{-3}$&$1.30532\times 10^{-2}$&$2.9694\times 10^{-3}$\\ \hline
M=3&$1.57535\times 10^{-3}$&$7.2083\times 10^{-3}$&$1.62779\times 10^{-3}$\\ \hline
M=5&$3.1167\times 10^{-4}$&$3.43748\times 10^{-3}$&$7.59478\times 10^{-4}$\\ \hline 
\end{tabular}
\end{center}
\caption{Relative errors of the Steklov approximations of $f_1, f_2,$ and $f_3$, respectively where $h=1$}
\end{table}

\begin{table}[here!]
\begin{center}
\begin{tabular}{|c||c|c|c| }
\hline
&rerr$_{2}(f_1)$&rerr$_{2}(f_2)$&rerr$_{2}(f_3)$\\ \hline
M=2&$5.13497\times 10^{-2}$&$1.69181\times 10^{-2}$&$2.77799\times 10^{-3}$\\ \hline
M=3&$4.15782\times 10^{-2}$&$1.0364\times 10^{-2}$&$1.52184\times 10^{-3}$\\ \hline
M=5&$1.78172\times 10^{-2}$&$4.58322\times 10^{-3}$&$6.98156\times 10^{-4}$\\ \hline 
\end{tabular}
\end{center}
\caption{Relative errors of the Steklov approximations of $f_1, f_2,$ and $f_3$, respectively where $h=0.8$}
\end{table}

\begin{table}[here!]
\begin{center}
\begin{tabular}{|c||c|c|c| }
\hline
&rerr$_{2}(f_1)$&rerr$_{2}(f_2)$&rerr$_{2}(f_3)$\\ \hline
M=2&$2.36676\times 10^{-1}$&$2.14194\times 10^{-2}$&$2.31158\times 10^{-3}$\\ \hline
M=3&$1.00467 \times 10^{-1}$&$1.04072\times 10^{-2}$&$1.3035\times 10^{-3}$\\ \hline
M=5&$5.79567\times 10^{-2}$&$5.45324\times 10^{-3}$&$5.9589\times 10^{-4}$\\ \hline 
\end{tabular}
\end{center}
\caption{Relative errors of the Steklov approximations of $f_1, f_2,$ and $f_3$, respectively where $h=0.5$}
\end{table}

It was observed that the above approximations were improved when some preliminary processing was performed. 
 In particular  it was worthwhile to first find the coefficients $a_j$ for a function $g_0(x,y) = a_0 + a_1 x+a_2y +a_3 xy$
 that interpolated the boundary data at the 4 corners of the rectangle.
 Then the Steklov appproximations  of solutions of Laplace's equation subject to the reduced boundary condition 
 $g_1(z) := g(z) - g_0(z)$ for $ z \in \bdy$ were observed to be better (have smaller error) than those for the boundary
 data $g$.
 
 Table (\ref{comp}) shows the comparison of relative errors of Steklov approximations of $f_1$ and $f_1+4$. Note that $f_1+4$ is the reduced boundary condition of $f_1$ such that  the value of the function at the 4 corners of $R_1$ is zero.

 \begin{table}[here!]
\begin{center}
\begin{tabular}{|c||c|c||c|c| }
\hline
&rerr$_{\infty}(f_1)$&rerr$_{\infty}(f_1+4)$&rerr$_{2}(f_1)$&rerr$_{2}(f_1+4)$\\ \hline
M=2&$6.59553\times 10^{-3}$&$5.27642\times 10^{-3}$&$5.22051\times 10^{-3}$&$2.54632\times 10^{-3}$\\ \hline
M=3&$2.28748\times 10^{-3}$&$1.82998\times 10^{-3}$&$1.57535\times 10^{-3}$&$7.6838\times 10^{-4}$\\ \hline
M=5&$5.55757\times 10^{-4}$&$4.46061\times 10^{-4}$&$3.1167\times 10^{-4}$&$1.52018\times 10^{-4}$\\ \hline 
\end{tabular}
\end{center}
\caption{Relative errors of Steklov approximations of $f_1$ and $f_1+4$ where $h=1$}\label{comp}

\end{table}

In his thesis Cho \cite{Ch} chapter 4, also investigated the approximation of harmonic functions by eigenfunctions of the Neumann
Laplacian on a rectangle. 
Even though such eigenfunctions form an orthogonal  basis of $\Hone$, finite approximations involving the first M eigenfunctions
were found to provide poor approximation properties for harmonic functions in $\Harm$.

\vspace{2em}
\section{ Approximations of Solutions of Robin Harmonic Boundary Value Problems } \label{s6}

When the first M harmonic Steklov eigenfunctions and eigenvalues are known, the associated Galerkin approximations
of Robin or Neumann boundary value problems for Laplace's equations may be found.
See Steinbach \cite{St}, chapter 8 or Zeidler \cite{Ze} chapter 19 for descriptions of such constructions and their general 
properties.
Here some specific error analyses  for harmonic functions will be proved and some numerical results will be described in 
the next section.

A function $u \in \Harm$ is said to be a (finite-energy) solution of the Robin harmonic boundary value problem on $\Om$
provided it satisfies
\beq \label{e6.1}
\IOm \gradu \cdot \gradv \dxy \pls b \ \Iby u \, v \dsg \eqs \Iby \, g \, v \dsg \foral v \in \Hone. \eeq

When (B1)  holds, standard variational arguments guarantee the existence and uniqueness of solutions 
of \eqref{e6.1} in $\Harm$. 
The solution is denoted $\Ebg$ and  satisfies the Robin boundary condition $ \Dnu u + bu \eqs g$ on $\bdy$ 
in a weak sense. 
For $b > 0$, it is 
\beq \label{e6.3}
\util(x,y) \eqs E_b g(x,y) \deqs \lim_{M \to \infty} \  \sum_{j=0}^M \ \frac{\gjh \sjt(x,y) }{b + \delj} \  \xfor (x,y) \in \Om. \eeq
This limit exists in the $H^1-$ norm provided $g \in H^{-1/2}(\bdy)$ as described in \cite{Au3}, section 10. 
In particular,  this holds when $g \in \Ltby$; note that even for linear functions on a rectangle, the Robin
or Neumann data $g$ may be discontinuous  on the boundary - so a useful analysis should allow such $g$.

When $g_M$ is given by \eqref{e3.11}, take $v = \sjt$ in \eqref{e6.1} to find that the solution is 
\beq \label{e6.5}
u_M(x,y) \deqs \Ebg_M \, (x,y) \eqs \frac{\gbar}{b} \pls \sum_{j=1}^M \ \frac{\gjh}{b + \deljh} \ \sjt(x,y) \xon \Omb. \eeq
Here $\deljh = \delj/|\bdy|$.
That is, after the Steklov spectrum has been found, the M-th Galerkin approximation of $E_b g$, just requires that the Steklov coefficients  $\gjh :=  \ang{g, \sjt}_{\bdy}$ be evaluated as in \eqref{e3.11}.

The error estimate for these approximations is the following \\
\btm \label{T6.1}
 Assume (B1) holds, $b > 0, \ g \in \Ltby$ and   $ g_M$ is  defined by \eqref{e3.11}.
 Then the function $u_M$ of \eqref{e6.5} is in $\Harm$ and 
 \beq \label{e6.6}
 \| \Ebg  - u_M \, \|_{\pal}^2 \leqs  \frac{1+\delta_{M+1}}{(b + \delta_{M+1})^2}\  \ \left[ \|g \|_{2, \bdy}^2  \mns \|g_M \|_{2, \bdy}^2 \right]
 \eeq  
 Moreover the functions $u_M$ converge uniformly to $\Ebg$ on compact subsets of $\Om$. \etm
\bpf
From \eqref{e6.1} and \eqref{e6.5}  one sees that 
\[ \Ebg(x,y) \mns \Ebg_M(x,y) \eqs \sum_{j = M+1}^{\infty} \ \frac{\gjh}{b + \deljh} \ \sjt(x,y) \xon \Omb \]

Evaluating the $\pal-$norm of this yields, using the orthogonality of the eigenfunctions, that
\[ \| \, \Ebg \mns \Ebg_M \, \|_{\pal}^2 \eqs  \sum_{j = M+1}^{\infty} \ \frac{\gjh^2 \ (1+ \delj)}{(b + \deljh)^2}  \] 

Thus
\beq \label {e6.7}
\| \, \Ebg \mns \Ebg_M \, \|_{\pal}^2 \leqs  \frac{ 1 + \delta_{M+1}}{(b+ \hat{\delta}_{M+1})^2} \quad  \| g - g_M \|_{2, \bdy}^2  \eeq
Since $\delta_M$ increase to infinity, the coefficient here is bounded so $\Ebg_M$ converges to $\Ebg$ in $\Hone$.
This equation implies \eqref{e6.6} as the Steklov eigenfunctions are $L^2-$ orthogonal on $\bdy$. 
Again the uniform convergence on compact subsets of $\Om$ is a standard result  for harmonic functions. 
\epf
 
 The estimate in \eqref{e6.6} shows again that $H^1-$error bounds for $\Ebg$ on $\Om$ may be found in terms of norms
 of $g-g_M$ on $\bdy$. 
 Some computational results for specific examples  are described in the next section.

 When the Neumann boundary condition ($b = 0$) holds then \eqref{e6.3} holds provided $\gbar = 0$ and  the solution 
is unique up to a constant.  
The minimum norm solution now is
\beq \label{e6.8}
\util(x,y) \eqs E_N g(x,y) \deqs \lim_{M \to \infty} \  \sum_{j=1}^M \ \frac{\gjh }{ \delj} \  \sjt(x,y) \  \xfor (x,y) \in \Om. \eeq
Let $u_M$  be this M-th partial sum, then $u_M$ converges to $E_N g$ in norm on $\Hone$ and $E_N$ is a continuous map of $H^{-1/2}(\bdy)$ to $\Harm$.
See section 10 of \cite{Au3} for more details.

The following error estimate for these approximations is proved using the same arguments as those for theorem \ref{T6.1}. 
\btm \label{T6.2}
 Assume (B1) holds, $ g \in \Ltby, \gbar = 0$ and   $ g_M$ is  defined by \eqref{e3.11}.
 Then  $u_M$ defined by  \eqref{e6.8} is in $\Harm$ and 
 \beq \label{e6.9}
 \| E_N g  - u_M \, \|_{\pal}^2 \leqs  \frac{1+\delta_{M+1}}{(b + \delta_{M+1})^2}\  \ \left[ \|g \|_{2, \bdy}^2  \mns \|g_M \|_{2, \bdy}^2 \right]
 \eeq  
 Moreover the functions $u_M$ converge uniformly to $\Ebg$ on compact subsets of $\Om$. \etm

\vspace{2em}
\section{Computation of solutions of Robin Harmonic Boundary value Problems} \label{s7}
The results of the preceding section provide representations of the solutions of Robin and Neumann problems for the 
Laplacian in terms of the harmonic Steklov eigenproblems.
Our observations are  that  approximations with relatively few (16-40) Steklov eigenfunctions compared quite well with
numerical solutions obtained using finite element software such as FreeFem++ (see \cite{He}).

Rather than compare the results with such software, however, we will present some data about comparisons with 
problems with exact solutions to illustrate the phenomenology observed.    
In particular we observed good approximations away from the boundary and some difficulty in handling
discontinuity in the data $g$ at points of discontinuity - even when the solution is nice. 
There is a Gibb's type effect in this case.

 Denote $\Gamo, \Gamtw, \Gamth$, and $\Gamf$ to be the side with $x=1, y=h, x=-1$, and $y=-h$, respectively such that $\partial\Omega= \Gamo \cup \Gamtw \cup \Gamth \cup \Gamf$.

\subsection{Neumann harmonic boundary value problem}~\\

Consider the boundary value problem on $\Om \eqs R_h$
\beq \label{nbvd}
\Delta \, u \eqs 0 \qquad \mbox{on $\R_h$ with} \quad  \Dnu \, u \eqs  \  g
\quad \mbox{on} \ \bdy
\eeq
with Dirichlet data 
\begin{equation}\label{bd1}
g(x,y)=
\begin{cases}
+1 \qquad &\mbox{on $\Gamo$ and $\Gamtw$}\\
-1 \qquad &\mbox{on $\Gamth$ and $\Gamf$}
\end{cases}
\end{equation}
We note that this example has a unique solution $u(x,y)=x+y$ with mean value zero on $R_h$. 
This solution is infinitely differentiable  but the boundary data $g$ is discontinuous at $(-1,h)$ and $(1,-h)$
because the domain $R_h$ has corners.

A graph of the numerical solution and of the error $u - u_5$ of the solution with $M=5$ is given in figure 2. 

\begin{table}[h!]
\begin{center}
\begin{tabular}{|c||c|c| }
\hline
&rerr$_{\infty}(u)$&rerr$_{2}(u)$\\ \hline
M=2&$3.44988\times 10^{-2}$&$2.17341\times 10^{-2}$\\ \hline
M=3&$2.34853\times 10^{-2}$&$1.23794\times 10^{-2}$\\ \hline
M=5&$1.43896\times 10^{-2}$&$5.98271\times 10^{-3}$\\ \hline 
\end{tabular}
\end{center}
\caption{Relative error of the Steklov approximation of the solution of (\ref{nbvd}) with the boundary condition (\ref{bd1}) where $h=1$}
\end{table}

\begin{figure}[h!]

\minipage{0.5\textwidth}
\includegraphics[height=9cm, width=\textwidth]{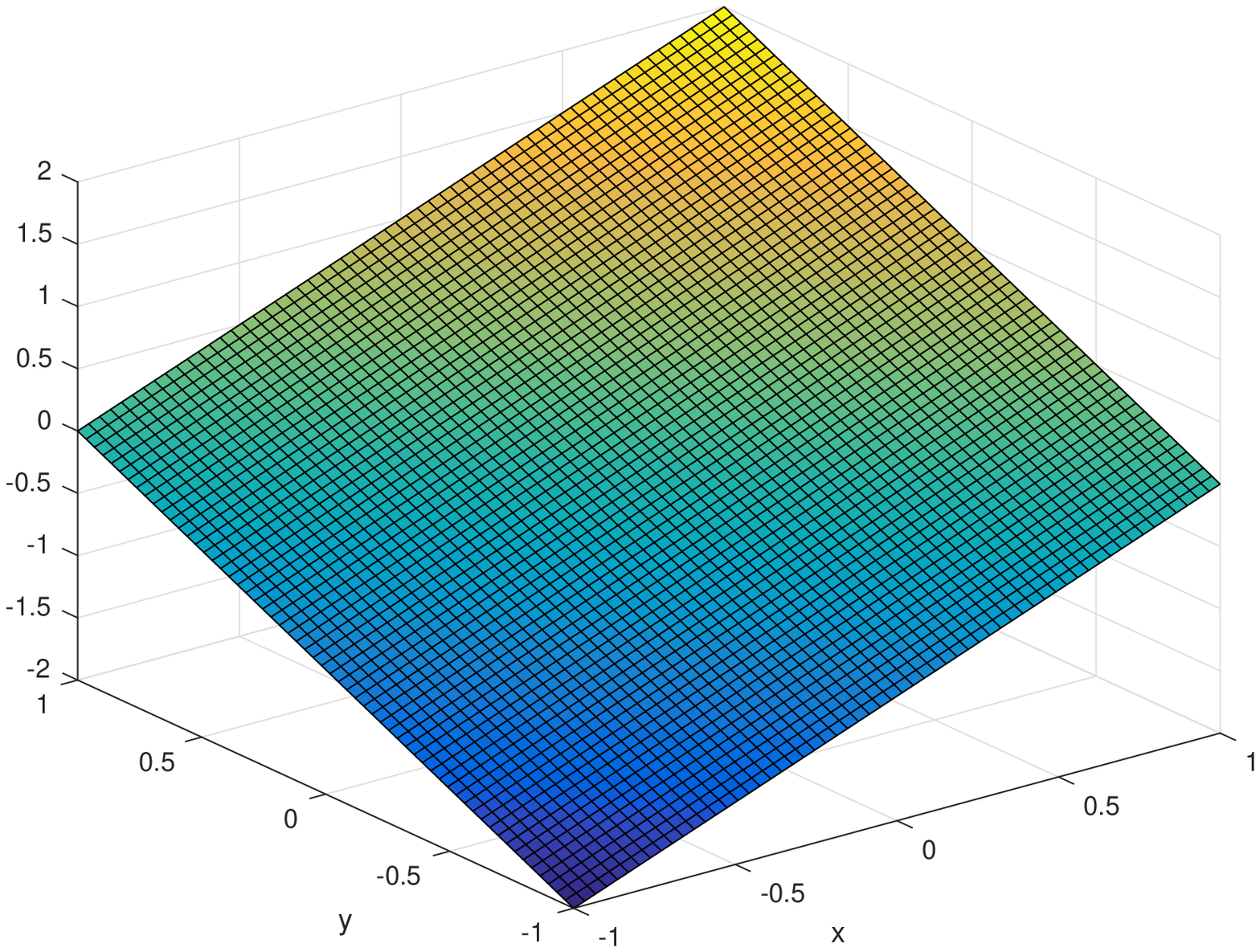}
\caption*{(a) Steklov approximation, $u_5$}
\endminipage
\minipage{0.5\textwidth}
\includegraphics[height=9cm, width=\textwidth]{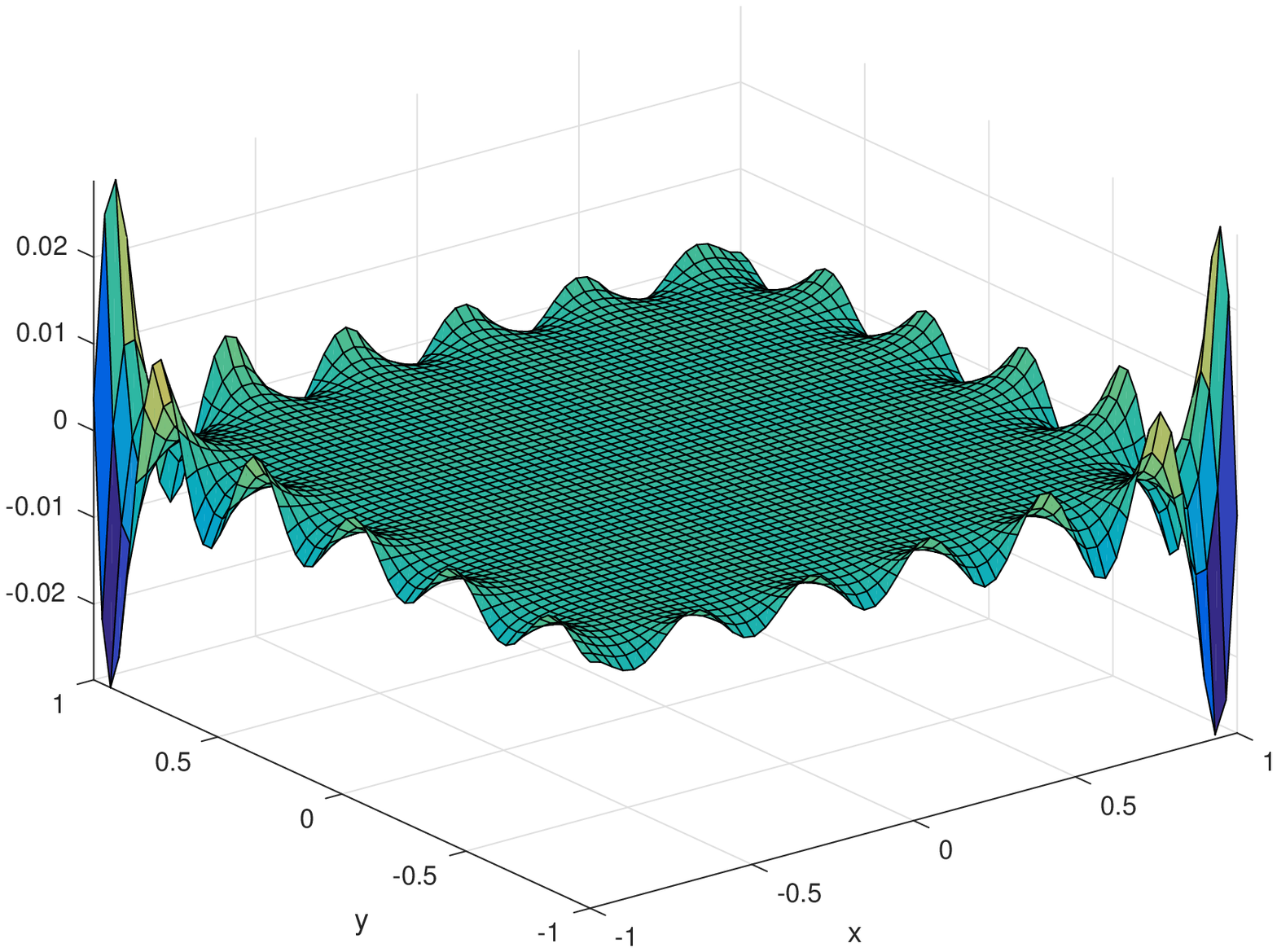}
\caption*{(b) Error in the solution, $u-u_5$}
\endminipage
\caption{Numerical results of the Steklov approximation of the solution of (\ref{nbvd}) with the boundary condition (\ref{bd1}) where $h=1$}
\label{lin}
\end{figure}

Another Neumann  problem (\ref{nbvd}) on $\Om \eqs R_h$ used  $g$ defined by
\begin{equation}\label{bd2}
g(x,y)=
\begin{cases}
+2 \qquad &\mbox{on $\Gamo$ and $\Gamth$}\\
-2h \qquad &\mbox{on $\Gamtw$ and $\Gamf$}
\end{cases}
\end{equation}

This problem has a unique solution $u(x,y)=x^2-y^2$ with mean value zero on the rectangle.
This solution is a well-known saddle function but now the boundary data $g$ is discontinuous at each corner.
Graphs of the Steklov approximation with $M = 5$ and the error $u - u_5$ are provided in figure 3.

\begin{table}[here!]
\begin{center}
\begin{tabular}{|c||c|c| }
\hline
&rerr$_{\infty}(u)$&rerr$_{2}(u)$\\ \hline
M=2&$9.07987\times 10^{-2}$&$1.32590\times 10^{-1}$\\ \hline
M=3&$5.34729\times 10^{-2}$&$9.20000\times 10^{-2}$\\ \hline
M=5&$2.64002\times 10^{-2}$&$5.70258\times 10^{-2}$\\ \hline 
\end{tabular}
\end{center}
\caption{Relative error of the Steklov approximation of the solution of (\ref{nbvd}) with the boundary condition (\ref{bd2}) where $h=1$}
\end{table}

\begin{figure}[h!]
\minipage{0.5\textwidth}
\includegraphics[height=9cm, width=\textwidth]{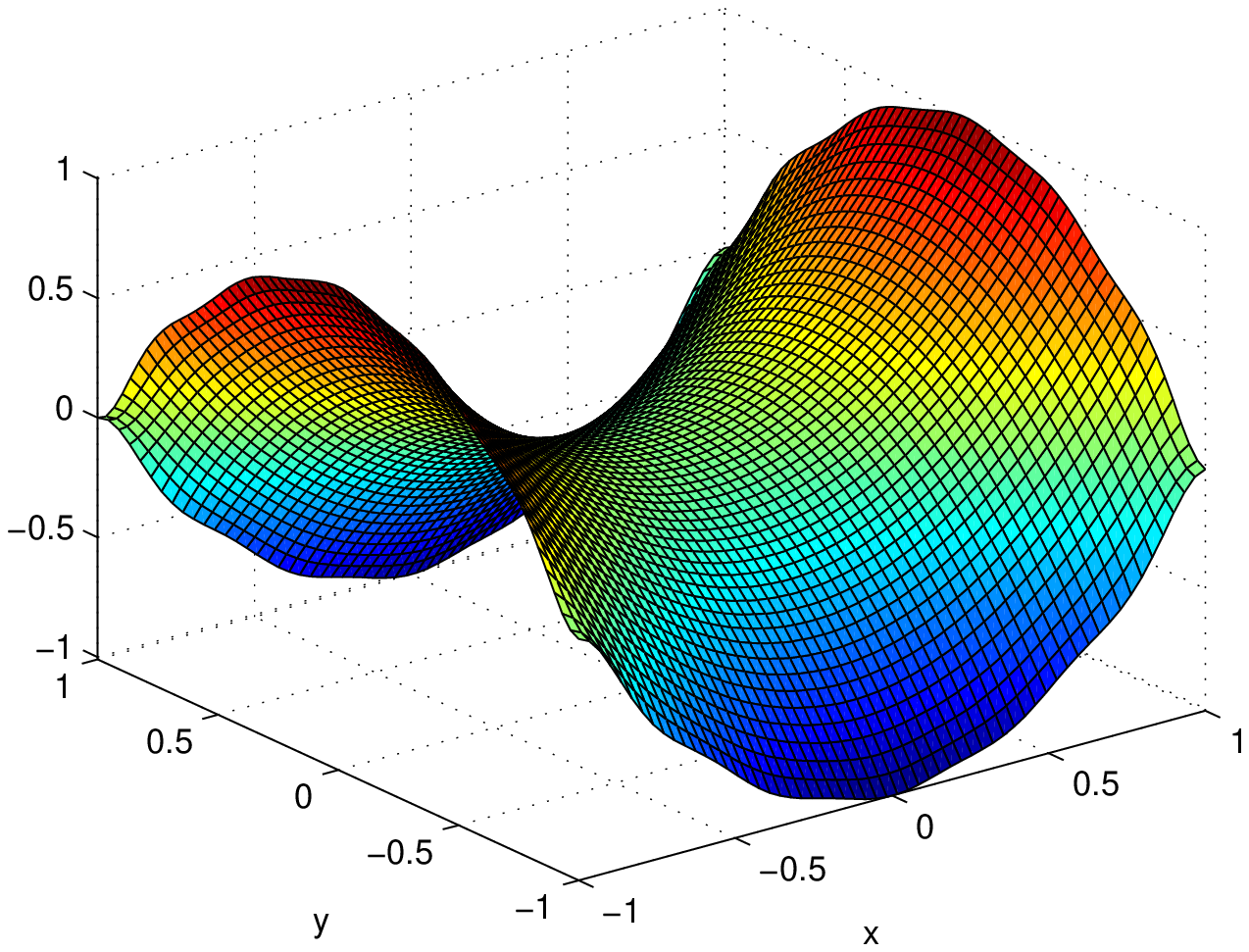}
\caption*{(a) Steklov approximation, $u_5$}
\endminipage
\minipage{0.5\textwidth}
\includegraphics[height=9cm, width=\textwidth]{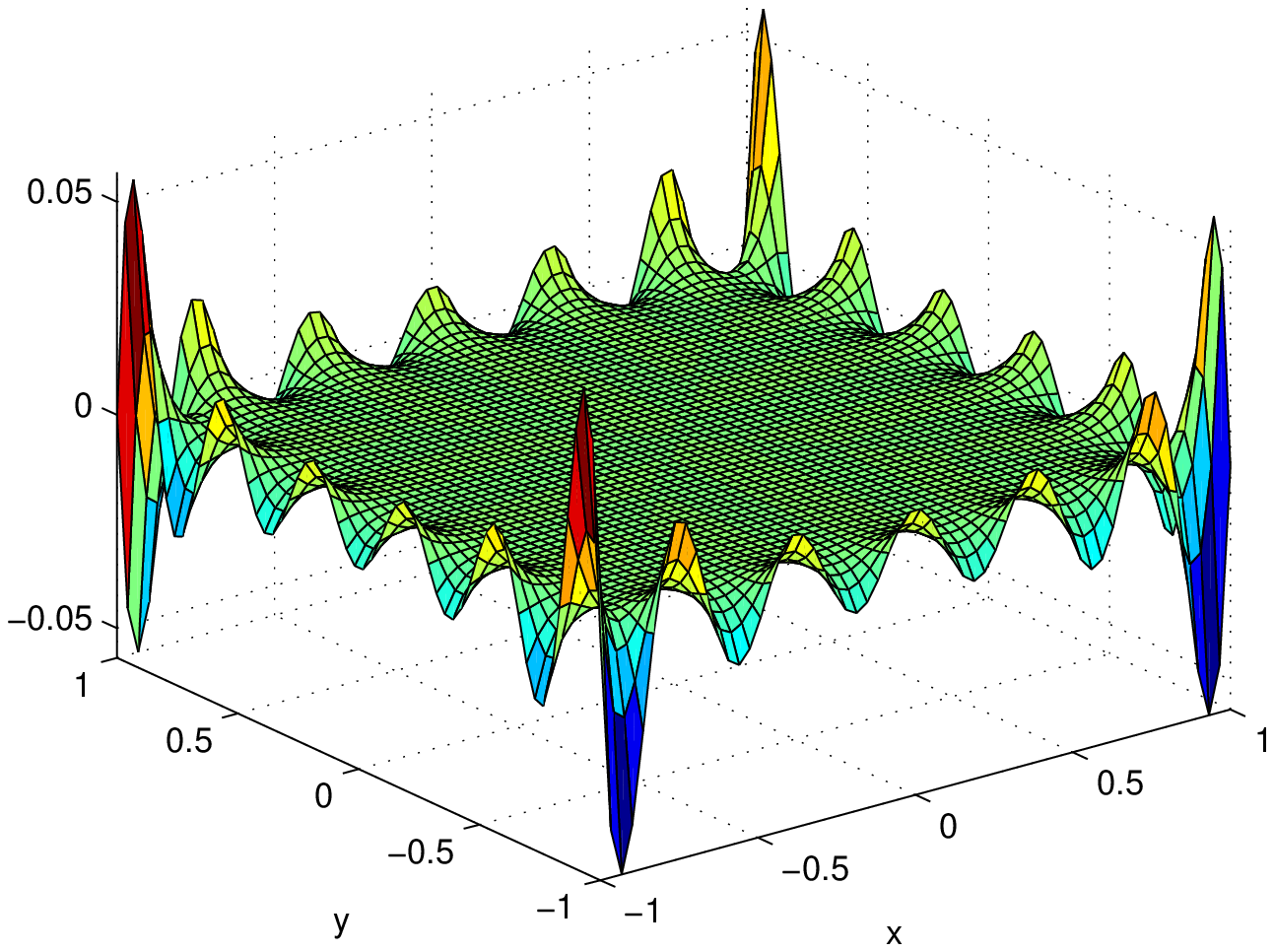}
\caption*{(b) Error in the solution, $u-u_5$}
\endminipage
\caption{Numerical results of the Steklov approximation of the solution of (\ref{nbvd}) with the boundary condition (\ref{bd2}) where $h=1$}
\label{lin}
\end{figure}

These simple examples show that the Steklov approximations of solutions of these problems provide
quite good approximations in the interior of the region even for small choices of M. 
The approximations satisfy the maximum principle, so the solutions are less accurate at, or near, the boundary.

\subsection{Robin harmonic boundary value problem}~\\

We consider a solution of the Robin harmonic boundary value problem with $b=1$ on $R_h$,
\beq\label{rbvd}
\Delta \, u \eqs 0 \qquad \mbox{on $R_h$ with} \quad   \Dnu u + bu \eqs g
\quad \mbox{on} \ \bdy
\eeq
 where $g$ is given by
\begin{equation}\label{bd3}
g(x,y)=
\begin{cases}
2(e^{1} \sin(y)) \qquad &\mbox{on $\Gamo$}\\
 e^x(\cos(h)+\sin(h)) \qquad &\mbox{on $\Gamtw$}\\
0 \qquad &\mbox{on $\Gamth$}\\
-e^x(\cos(h)+\sin(h)) \qquad &\mbox{on $\Gamf$}
\end{cases}
\end{equation}

The unique solution of this problem is $u(x,y)=e^{x}\sin(y)$.
The Steklov approximation with $M = 5$ is shown in figure 4,  together with a graph of the error function
$u - u_5$. 
Again the relative error is quite reasonable and the approximations are very accurate away from the 
boundary.

\begin{table}[here!]
\begin{center}
\begin{tabular}{|c||c|c| }
\hline
&rerr$_{\infty}(u)$&rerr$_{2}(u)$\\ \hline
M=2&$1.51186\times 10^{-2}$&$1.4854\times 10^{-2}$\\ \hline
M=3&$9.64123 \times 10^{-3}$&$7.84911\times 10^{-3}$\\ \hline
M=5&$5.60122\times 10^{-3}$&$3.53263\times 10^{-3}$\\ \hline 
\end{tabular}
\end{center}
\caption{Relative error of the Steklov approximation of the solution of (\ref{rbvd}) with the boundary condition (\ref{bd3}) where $h=1$}
\end{table}

\begin{figure}[h!]
\minipage{0.5\textwidth}
\includegraphics[height=9cm, width=\textwidth]{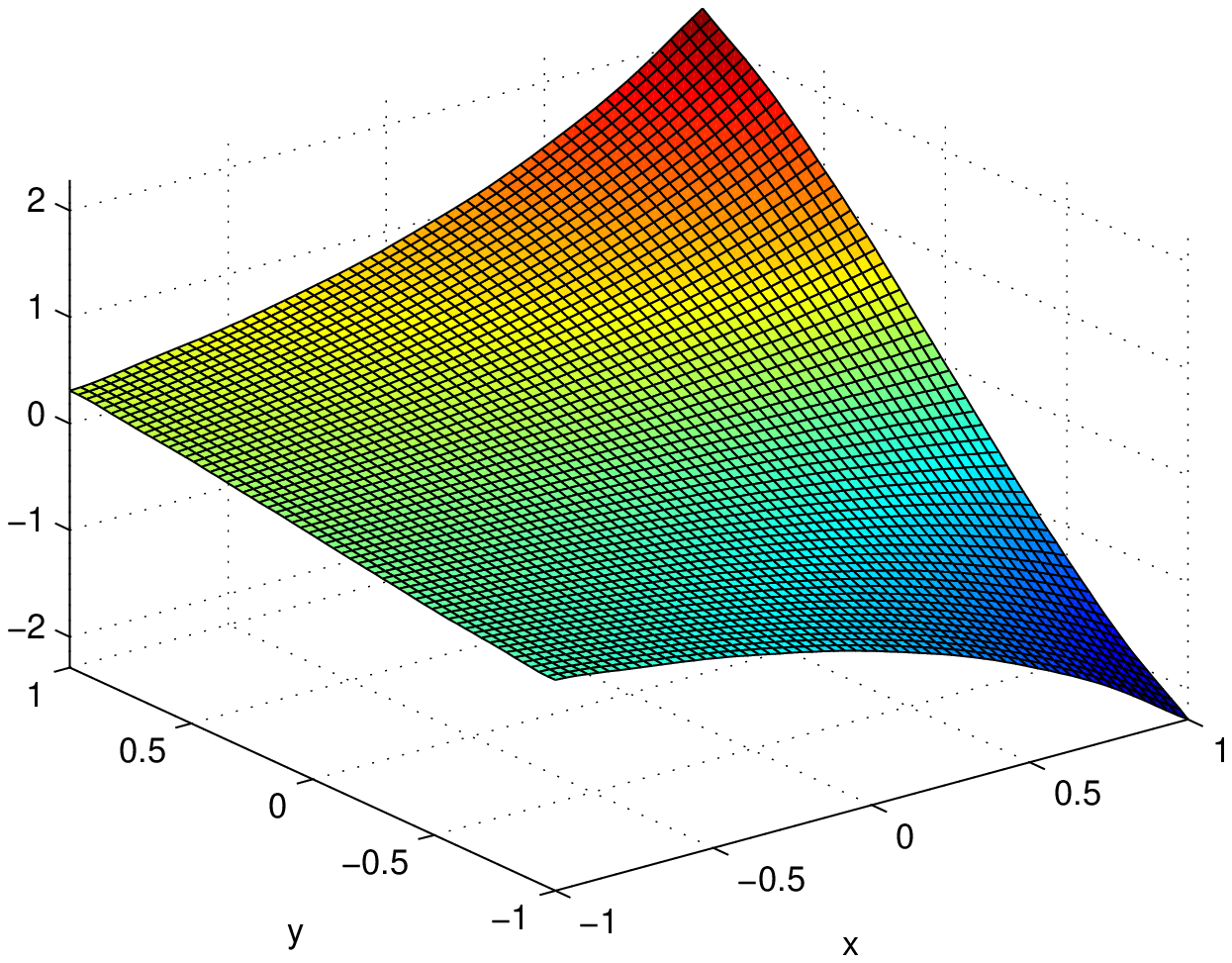}
\caption*{(a) Steklov approximation, $u_5$}
\endminipage
\minipage{0.5\textwidth}
\includegraphics[height=9cm, width=\textwidth]{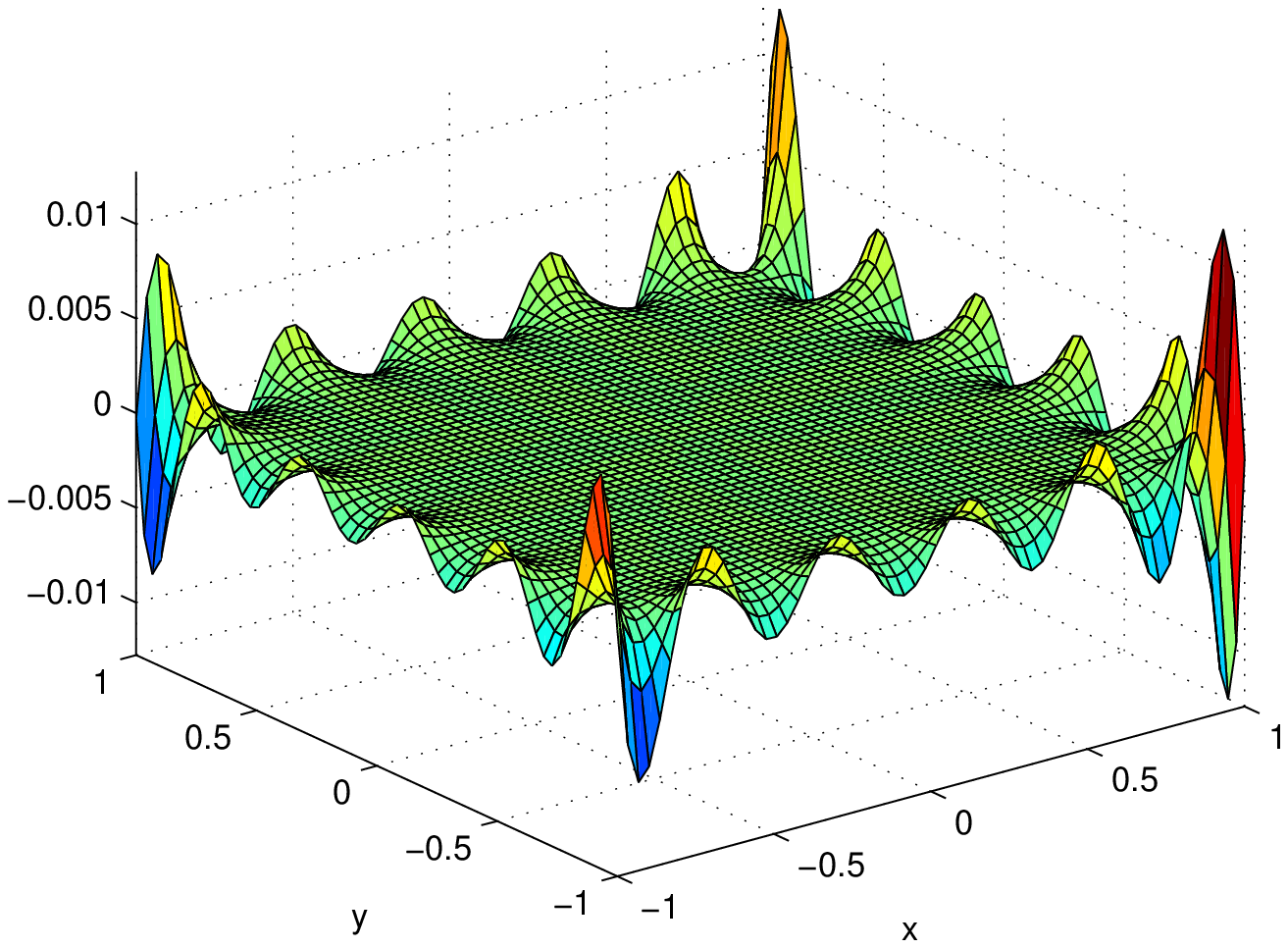}
\caption*{(b) Error in the solution, $u-u_5$}
\endminipage
\caption{Numerical results of the Steklov approximation of the solution of (\ref{rbvd}) with the boundary condition (\ref{bd3}) where $h=1$}
\label{lin}
\end{figure}

These simple examples were chosen primarily to illustrate the phenomenology observed in computing Steklov approximations.
There clearly are many further questions about the efficacy of such approximations but the primary observation is 
that low order Steklov approximations do provide good interior approximations to solutions of harmonic boundary value problems.

%
%

\end{document}